\documentclass{amsart}
\DeclareMathSymbol{\twoheadrightarrow} 
{\mathrel}{AMSa}{"10}

\def\Q{{\mathbf Q}}
\def\Z{{\mathbf Z}}

\def\R{{\mathfrak  R}}
\def\F{{\mathbf F}}
\def\M{{\mathbf M}}
\def\L{{\mathbf L}}

\def\Sz{\mathbf{Sz}}

\def\Sn{{\mathbf S}_n}
\def\An{{\mathbf A}_n}
\def\bchi{{\mathbf \chi}}
\def\bphi{{\mathbf \phi}}
\def\Gal{\mathrm{Gal}}
\def\Perm{\mathrm{Perm}}

\def\End{\mathrm{End}}
\def\Aut{\mathrm{Aut}}

\def\I{\mathrm{Id}}

\def\fchar{\mathrm{char}}

\def\GL{\mathrm{GL}}

\def\Sp{\mathrm{Sp}}
\def\PSL{\mathrm{PSL}}
\def\PGL{\mathrm{PGL}}
\def\M{\mathrm{M}}
\def\dim{\mathrm{dim}}

\def\P{{\mathbf P}}

\newtheorem{thm}{Theorem}[section]
\newtheorem{lem}[thm]{Lemma}

\newtheorem{prop}[thm]{Proposition}
\theoremstyle{definition}
\newtheorem{defn}[thm]{Definition}
\newtheorem{ex}[thm]{Example}

\newtheorem{rem}[thm]{Remark}
\newtheorem{rems}[thm]{Remarks}
\title{Hyperelliptic jacobians and  projective linear Galois groups}
\author[Yuri G. Zarhin]{Yuri G. Zarhin}
\address{Department of Mathematics, Pennsylvania State University,
University Park, PA 16802, USA}
\email{zarhin\char`\@math.psu.edu}
\thanks{Partially supported by NSF grant DMS-0070664}

\begin{document}
\maketitle
\section{Introduction}

In \cite{Zarhin} the author  proved that
in characteristic $0$ the
jacobian $J(C)=J(C_f)$  of a hyperelliptic curve
$$C=C_f:y^2=f(x)$$
has only trivial endomorphisms
over an algebraic closure $K_a$ of the ground field $K$
if the Galois group $\Gal(f)$ of the irreducible polynomial
$f \in K[x]$ is ``very big". Namely, if $n=\deg(f) \ge 5$
and $\Gal(f)$ is either the symmetric group $\Sn$ or the alternating group $\An$
 then the ring $\End(J(C_f))$ of $K_a$-endomorphisms of $J(C_f)$ coincides with $\Z$.  Later the author \cite{Zarhin2} proved that $\End(J(C_f))=\Z$ for an infinite series 
of $\Gal(f)=\L_2(2^{r}):=\PSL_2(\F_{2^r})$ and $n=2^{r}+1$ 
(with $r \ge 3$ and $\dim(J(C_f))=2^{r-1}$) or when $\Gal(f)$ is 
 the Suzuki group $\Sz(2^{2r+1})$ and $n=2^{2(2r+1)}+1$ 
(with $\dim(J(C_f))=2^{4r+1}$). He also proved the same assertion when $n=11$ or $12$ and $\Gal(f)$ is the Mathieu group $\M_{11}$ or $\M_{12}$.
 (In those cases $J(C_f)$ has dimension $5$.)

We refer the reader to \cite{Mori1}, \cite{Mori2}, 
\cite{Katz1}, \cite{Katz2}, \cite{Masser}, \cite{Zarhin}, \cite{Zarhin2} 
for a discussion of known results about, and examples  
of, hyperelliptic jacobians without complex multiplication.

In the present paper we prove that 
$\End(J(C_f))=\Z$ when the set $\R=\R_f$ of roots of $f$ can be identified with
the $(m-1)$-dimensional projective space $\P^{m-1}(\F_q)$ over a finite field $\F_q$ of {\sl odd} characteristic in such a way that $\Gal(f)$, viewed as a permutation group of $\R_f$, becomes either the projective linear group $\PGL(m,\F_q)$ or the projective special linear group $\L_m(q):=\PSL(m,\F_q)$. Here we assume that $m>2$. In this case
$$n=\deg(f)=\#(\P^{m-1}(\F_q))=\frac{q^{m}-1}{q-1}$$
 and 
$\dim(J(C_f))$ is $[(\frac{q^{m}-1}{q-1}-1)/2]$, i.e. the integral part of  $(\frac{q^{m}-1}{q-1}-1)/2$.

Our proof is based on a result of Guralnick \cite{GurTiep}, who proved that in the ``generic" case the dimension of each nontrivial irreducible representation of $\L_m(q)$ in characteristic $2$ is greater than or equal to 
$$2[(\frac{q^{m}-1}{q-1}-1)/2].$$

We also discuss the similar problem when $K$ has prime chatacteristic $>2$. It turns out that $\End(J(C_f))=\Z$ under an additional assumption that $m$ is {\sl even} (i.e., when $n$ is even).
The case of $n=12$ and $\Gal(f)=\M_{12}$ is also treated.

\section{Main results}
\label{mainr}
Throughout this paper we assume that $K$ is a field with $\fchar(K) \ne 2$. We fix its algebraic closure $K_a$ and write $\Gal(K)$ for the absolute Galois group $\Aut(K_a/K)$. If $X$ is an abelian variety defined over $K$ then we write $\End(X)$ for the ring of  $K_a$-endomorphisms of $X$.

Suppose $f(x) \in K[x]$ is a separable polynomial 
 of degree $n\ge 5$. Let $\R=\R_f \subset K_a$ be 
 the set of roots of $f$, let $K(\R_f)=K(\R)$ be 
 the splitting field of $f$ and $\Gal(f):=\Gal(K(\R)/K)$ the Galois group of $f$, viewed as a subgroup of $\Perm(\R)$.
Let $C_f$ be the hyperelliptic curve $y^2=f(x)$. Let  $J(C_f)$ be
its jacobian, $\End(J(C_f))$ the ring of $K_a$-endomorphisms of $J(C_f)$.

\begin{thm}
\label{main2}
Assume that there exist a positive integer $m>2$ and an odd power prime $q$ such that $n=\frac{q^m-1}{q-1}$ and $\Gal(f)$ contains a subgroup isomorphic  to $\L_m(q)$. (E.g., $\Gal(f)$ is isomorphic either to
$\PGL_m(\F_q)$ or to $\L_m(q)=\PSL_m(\F_q)$.)

Then either $\End(J(C_f))=\Z$ or $m$ is odd, $\fchar(K)>0$ and
$J(C_f)$ is a supersingular abelian variety.
\end{thm}

\begin{rem}
\label{parity}
Clearly $m$ is even if and only if $n$ is even.
\end{rem}

\begin{rem}
\label{enlarge}
Replacing $K$ by $K(\R)^{\L_m(q)}$, we may assume that
$$\Gal(f)=\L_m(q).$$
Also, taking into account that $\L_m(q)$ is simple non-abelian and
replacing $K$ by its abelian extension obtained by adjoining all $2$-power roots of unity, we may assume that $K$ contains all $2$-power roots of unity.
\end{rem}

\begin{thm}
\label{M12}
Suppose $n=12$ and $\Gal(f)$ is isomorphic to the Mathieu group $\M_{12}$.
Then $\End(J(C_f))=\Z$.
\end{thm}

\begin{rem}
\label{enlargeM}
When $\fchar(K)=0$ the assertion of Theorem \ref{M12} is proven in \cite{Zarhin2}. Taking into account that $\M_{12}$ is simple non-abelian and
replacing $K$ by its abelian extension obtained by adjoining all $2$-power roots of unity, we may assume that $K$ contains all $2$-power roots of unity.
\end{rem}

We will prove Theorems \ref{main2} and \ref{M12} in \S \ref{final}.

\section{Permutation groups and permutation modules}
\label{permute}

Let $B$ be a finite set consisting of $n \ge 5$ elements. We write $\Perm(B)$ for the group of permutations of $B$. A choice of ordering on $B$ gives rise to an isomorphism
$$\Perm(B) \cong \Sn.$$
Let $G$ be a  subgroup of $\Perm(B)$.
For each $b \in B$ we write $G_b$ for the stabilizer of $b$ in $G$; it is a subgroup of $G$.

\begin{rem}
\label{transitive}
Assume that the action of $G$ on $B$ is transitive.
It is well-known that each $G_b$ is  of index $n$ in $G$ and 
 all the $G_b$'s are conjugate in $G$. 
 Each conjugate of $G_b$ in $G$ is the stabilizer of a point of $B$. 
 In addition, one may identify the $G$-set $B$ with the set of cosets $G/G_b$ 
with the standard action by $G$. 
\end{rem}

Let $\F$ be a field. We write $\F^B$
for the $n$-dimensional $\F$-vector space of maps $h:B \to \F$.
The space $\F^B$ is  provided with a natural action of $\Perm(B)$ defined
as follows. Each $s \in \Perm(B)$ sends a map
 $h:B\to \F$ into  $sh:b \mapsto h(s^{-1}(b))$. The permutation module $\F^B$
contains the $\Perm(B)$-stable hyperplane
$$(\F^B)^0=
\{h:B\to\F\mid\sum_{b\in B}h(b)=0\}$$
and the $\Perm(B)$-invariant line $\F \cdot 1_B$ where $1_B$ is the constant function $1$. The quotient $\F^B/(\F^B)^0$ is a trivial $1$-dimensional $\Perm(B)$-module.

Clearly, $(\F^B)^0$ contains $\F \cdot 1_B$ if and only if $\fchar(\F)$ divides $n$. If this is {\sl not} the case then there is a $\Perm(B)$-invariant splitting
$$\F^B=(\F^B)^0 \oplus \F \cdot 1_B.$$

Clearly, $\F^B$ and $(\F^B)^0$  carry natural structures of $G$-modules. Their characters depend only on the characteristic of $\F$.

Let us consider the case of $\F=\Q$. Then the character of $\Q^B$
 is called the {\sl permutation character} of $B$. Let us denote by 
$\bchi=\bchi_B:G \to \Q$ the character of $(\Q^B)^0$.
 Clearly, $1+\bchi$ is the permutation character of $B$.

Now, let us consider the case of $\F=\F_2$.
If  $n$ is even then let us define the $\Perm(B)$-module 
$$Q_B:=(\F_2^B)^0/(\F_2 \cdot 1_B).$$
If $n$ is odd then let us put
$$Q_B:=(\F_2^B)^0.$$

\begin{rem}
Clearly, $Q_B$ is a faithful $G$-module. If $n$ is odd then $\dim_{\F_2}(Q_B)=n-1$. If  $n$ is even then $\dim_{\F_2}(Q_B)=n-2$.
\end{rem}

 Let $G^{(2)}$ be the set of $2$-regular elements of $G$. 
 Clearly, the Brauer character of the $G$-module $\F_{2}^B$  
 coincides with the restriction of $1+\bchi_B$ to $G^{(2)}$. 
 This implies easily that the Brauer character of the $G$-module $(\F_{2}^B)^0$ 
 coincides with the restriction of   $\bchi_B$ to $G^{(2)}$.

\begin{rem}
\label{Bcharacter}
 Let us denote by 
$\bphi_B=\bphi$
 the Brauer character of the $G$-module $Q_B$.
 One may easily check that $\bphi_B$ coincides with the restriction of 
 $\bchi_B$ to $G^{(2)}$ if  $n$ is odd and with the restriction of 
$\bchi_B-1$ to $G^{(2)}$ if  $n$ is even. 
\end{rem}

We refer to \cite{Zarhin2} for a discussion of the following definition.

\begin{defn}
Let $V$ be a vector space over a field $\F$, let $G$ be a group and
$\rho: G \to \Aut_{\F}(V)$ a linear representation of $G$ in $V$. We
say that the $G$-module $V$ is {\sl very simple} if it enjoys the
following property:

If $R \subset \End_{\F}(V)$ is an $\F$-subalgebra containing the
identity operator $\I$ such that

 $$\rho(\sigma) R \rho(\sigma)^{-1} \subset R \quad \forall \sigma \in G$$
 then either $R=\F\cdot \I$ or $R=\End_{\F}(V)$.
\end{defn}

\begin{rems}
\label{image}
\begin{enumerate}
\item[(i)]
If $G'$ is a subgroup of $G$ and the $G'$-module $V$ is very simple 
then obviously the $G$-module $V$ is also very simple.

\item[(ii)]
A  very simple module is absolutely simple (see \cite{Zarhin2}, Remark 2.2(ii)).

\item[(iii)]
If $\dim_{\F}(V)=1$ then obviously the $G$-module $V$ is very simple.

\item[(iv)]
Assume that the $G$-module $V$ is  very simple and $\dim_{\F}(V)>1$.
Then  $V$ is not induced from a subgroup $G$ (except $G$ itself) and is not isomorphic to a tensor product of two $G$-modules, whose $\F$-dimension is strictly less than $\dim_{\F}(V)$  (see \cite{Zarhin2}, Examples 7.1).

\item[(v)] If $\F=\F_2$ and $G$ is {\sl perfect} then  properties (ii)-(iv)
characterize the very simple $G$-modules (see \cite{Zarhin2}, Th. 7.7).
\end{enumerate}
\end{rems}

The following statement provides a criterion of very simplicity over $\F_2$.

\begin{thm}
\label{Very3}
Suppose  a positive integer $N>1$ and a group $H$ enjoy the following properties:
\begin{itemize}
\item
$H$ does not contain a subgroup of index dividing $N$ except $H$ itself.

\item
Let $N=ab$ be a factorization of $N$ into a product of two
positive integers $a>1$ and $b>1$. Then either
there does not exist an absolutely simple $\F_2[H]$-module of $\F_2$-dimension $a$ or
there does not exist an absolutely simple $\F_2[H]$-module of $\F_2$-dimension $b$.
\end{itemize}

Then each absolutely simple $\F_2[H]$-module of $\F_2$-dimension $N$ is very simple. 
\end{thm}

\begin{proof}
This is Corollary 4.12 of \cite{Zarhin2}.
\end{proof}

\begin{thm}
\label{Lnq}
Suppose that there exist a positive integer $m>2$ and an odd power prime $q$ such that $n=\frac{q^m-1}{q-1}$. Suppose $G$ is a subgroup of $\Sn$. Suppose $G$ contains a subgroup isomorphic  to $\L_m(q)$. Then the $G$-module $Q_B$ is very simple.
\end{thm}

The rest of this section is devoted to the proof of Theorem \ref{Lnq}.
 In light of Remark \ref{image}(ii), we may assume that $G=\L_m(q)$.

{\bf ``Generic" case}.
Assume that $(m,q) \ne (4,3)$. Then it follows from Theorem 1.1 (applied to $p=2$) and the Table III of \cite{GurTiep} that each nontrivial (absolutely) irreducible representation of $\L_m(q)$ in characteristic $2$ has dimension which is greater or equal than $N:=\dim_{\F_2}(Q_B)$. Taking into account that $\L_m(q)$ is (simple) not solvable and $Q_B$ is a faithful $\L_m(q)$-module, we conclude that $Q_B$ is absolutely simple. 

Now we claim that the group $G=\L_m(q)$ does not contain a subgroup of index dividing $N:=\dim_{\F_2}(Q_B)$ except $G$ itself.

Indeed, if $G'$ is a subgroup of $G$ such that $G' \ne G$ and $[G:G']$ divides $\dim_{\F_2}(Q_B)$ then the simple group $G$ acts faithfully on $B'=G/G'$ and therefore $[G:G'] \ge 5$. In particular, we get a faithful $G$-module $Q_{B'}$, whose dimension is strictly less than $\dim_{\F_2}(Q_B)$.  

Since each strict divisor $a$ of $N$ lies strictly between $1$ and $N$,
there does not exist an absolutely simple $\F_2[G]$-module of $\F_2$-dimension $a$.

Now the very simplicity of the $G$-module $Q_B$ follows from Theorem \ref{Very3}.

{\bf The special case of} $m=4,q=3$. We have $n=\#(B)=40$ and $\dim_{\F_2}(Q_B)=38$. 
According to the Atlas (\cite{Atlas}, pp. 68-69), $G=\L_4(3)$
has two conjugacy classes of maximal subgroups  of index $40$. 
All other maximal subgroups have index greater than $40$.
 Therefore all subgroups of $G$ (except $G$ itself) have index greater than $39>38$.
This implies that each action of $G$ on $B$ is transitive.
The permutation character (in notations of \cite{Atlas}) is (in both cases) $1+\chi_4$, i.e., $\bchi=\chi_4$.
Since $40$ is even, we need to consider 
the restriction of $\bchi-1$ to the set of $2$-regular elements of $G$ 
and this restriction coincides with the absolutely irreducible Brauer character
 $\phi_4$ (in notations of \cite{AtlasB}, p. 165).
 In particular, the corresponding $G$-module $Q_B$ is absolutely simple.
It follows from the Table on p. 165 of \cite{AtlasB} 
that all absolutely irreducible representations
 of $G$ in characteristic $2$ have dimension which is {\sl not}
 a strict divisor of $38$. 
Combining this observation with the absence of subgroups in $G$ 
of index less or equal than $38$, we conclude, 
thanks to Theorem \ref{Very3}, that $Q_B$ is very simple.
This ends the proof of Theorem \ref{Lnq}.

\section{Proof of Theorems \ref{main2} and \ref{M12}}
\label{final}
Recall that $\Gal(f) \subset \Perm(\R)$. In addition, it is known that the natural homomorphism $\Gal(K) \to \Aut_{\F_2}(J(C)_2)$ factors 
through the canonical surjection $\Gal(K) \twoheadrightarrow \Gal(K(\R)/K)=\Gal(f)$
and the $\Gal(f)$-modules $J(C)_2$ and $Q_{\R}$ are isomorphic (see, for instance, Th. 5.1 of \cite{Zarhin2}). In particular, if the $\Gal(f)$-module $Q_{\R}$ is very simple then the $\Gal(f)$-modules $J(C)_2$ is also very simple and therefore is absolutely simple.

\begin{lem}
\label{cor51}
If the $\Gal(f)$-module $Q_{\R}$ is very simple then either $\End(J(C_f))=\Z$ or $\fchar(K)>0$ and $J(C_f)$ is a supersingular abelian variety.
\end{lem}

\begin{proof}
This is Corollary 5.3 of \cite{Zarhin2}.
\end{proof}

It follows from Theorem \ref{Lnq} that under the assumptions of Theorem \ref{main2}, the $\Gal(f)$-module $Q_{\R}$ is very simple. Applying Lemma \ref{cor51}, we conclude that  either $\End(J(C_f))=\Z$ or $\fchar(K)>0$ and $J(C_f)$ is a supersingular abelian variety. 

If $n=12$ and $\Gal(f) \cong \M_{12}$ then the $\Gal(f)$-module $Q_B$ is also very simple (\cite{Zarhin2}, Th. 7.12(ii)). Again we conclude that under the assumptions of Theorem  \ref{M12} either $\End(J(C_f))=\Z$ or $\fchar(K)>0$ and $J(C_f)$ is a supersingular abelian variety (\cite{Zarhin2}, Th. 7.13(ii)).

In order to finish the proof  of Theorem \ref{main2}  we need only to check that $J(C_f))$ is {\sl not} supersingular if $m$ is even. Similarly, in order to prove Theorem \ref{M12} we need only
to check that if  $(n,\Gal(f))=(12, \M_{12})$ then
$J(C_f))$ is {\sl not} supersingular.
 Using Remarks \ref{enlarge} and \ref{enlargeM}, we may assume that either $\Gal(f)=\L_m(q)$ or $(n,\Gal(f))=(12, \M_{12})$ and in both cases $K$ contains all $2$-power roots of unity.
Clearly, the desired assertions are immediate corollaries of the following statement.

 \begin{lem}
\label{supernot}
Suppose an even positive integer $n$ and a finite simple non-abelian group $G$ enjoy one of the following two properties.
\begin{enumerate}
\item[(i)]
There exist  an odd power prime $q$ and an even integer $m \ge 4$ such that
$n=(q^m-1)/(q-1)$ and $G \cong \L_m(q)$;
\item[(ii)]
$n=12$ and $G \cong \M_{12}$.
\end{enumerate}

Let us put $g=(n-2)/2$.
Suppose $F$ is a field, whose characteristic is not $2$. Suppose that $F$  contains all $2$-power roots of unity. Suppose that $X$ is a $g$-dimensional abelian variety  over  $F$ such that the image of $\Gal(F)$ in $\Aut(X_2)$ is isomorphic to $G$ and the  $G$-module $X_2$ is absolutely simple. Then $X$ is not supersingular.
\end{lem}

\begin{proof}[Proof of Lemma \ref{supernot}]
  Every nontrivial representation of $G$ in characteristic $2$ has dimension $>g$. Indeed, first assume that $G=\L_m(q)$. Then
in the ``generic" case" of $(m,q) \ne (4,3)$ such a representation must have dimension $\ge 2g>g$, thanks to the already cited Th. 1.1 and Table III of \cite{GurTiep}. If $(m,q) = (4,3)$ then $n=40, 2g=38$ and the smallest dimension is $26>19=g$, according to the Tables in \cite{AtlasB}. Second, if $G=\M_{12}$ then this assertion  follows from Th. 8.1 on p. 80 in  \cite{James}; see also the Tables in \cite{AtlasB}. 

\begin{prop}
\label{double}
 Suppose $G'\twoheadrightarrow G$ is a central extension of $G$. In addition, assume that either $G'=G$ or $G'$ is a double cover of $G$, i.e., $\ker(G' \twoheadrightarrow G)$ is a central cyclic subgroup of order $2$ in $G'$. Suppose $V$ is a finite-dimensional $\Q_2$-vector space and
$$\rho: G'\to\Aut_{\Q_2}(V)$$ 
is an absolutely irreducible faithful representation of $G'$ over $\Q_2$. Then
	$$\dim_{\Q_2}(V) \ne 2g.$$
\end{prop}

\begin{proof}[Proof of Proposition \ref{double}]
Clearly, $\rho$ defines an absolutely irreducible projective representation of $G$ in $V$ over $\Q_2$. 

Assume first that $G=\L_m(q)$. Then in the ``generic" case 
every absolutely irreducible nontrivial projective representation of $G$ in characteristic $0$ must have dimension $\ge 2g+1>2g$ (see \cite{GurTiep}, Table II). If $(m,q) = (4,3)$ then the Proposition follows from the Tables in \cite{Atlas}.

Second, suppose $G=\M_{12}$. Then $n=12, 2g=10$. All faithful absolutely irreducible  representations of $\M_{12}$ in characteristic zero have dimension $\ge 11>10$ (\cite{Atlas}, p. 33). This  proves the Proposition in the case when $G'=G=\M_{12}$ and also when $G'$ is a trivial double cover, i.e., is isomorphic to a product of $\M_{12}$ and a cyclic group of order $2$. If $G'$ is a  nontrivial double cover of $\M_{12}$ then it has precisely two non-isomorphisc $10$-dimensional absolutely irreducible representations in characteristic $0$ (up to an isomorphism) \cite{Hump}. However, none of them is defined over $\Q_2$. Indeed, each character of $G'$ of degree $10$  takes on a value, whose square is $-2$  (\cite{Hump}, Table 1 on p. 410; \cite{Atlas}, p. 33).
\end{proof}

Assume that $X$ is supersingular. Our goal is to get a contradiction.
We write $T_2(X)$ for the $2$-adic Tate module of $X$
and
$$\rho_{2,X}:\Gal(F) \to \Aut_{\Z_2}(T_2(X))$$
for 
the corresponding $2$-adic representation. It is well-known that $T_2(X)$ is
a free $\Z_2$-module of rank $2\dim(X)=2g$ and
$$X_2=T_2(X)/2 T_2(X)$$
(as Galois modules). Let us put
$$H=\rho_{2,X}(\Gal(F)) \subset \Aut_{\Z_2}(T_2(X)).$$
Clearly, the natural homomorphism
$$\bar{\rho}_{2,X}:\Gal(F) \to \Aut(X_2)$$
defining the Galois action on the points of order $2$ is the composition of
$\rho_{2,X}$ and the (surjective) reduction map modulo $2$
$$\Aut_{\Z_2}(T_2(X)) \twoheadrightarrow \Aut(X_2).$$
This gives us a natural (continuous) {\sl surjection}
$$\pi:H \twoheadrightarrow \bar{\rho}_{2,X}(\Gal(F)) \cong G,$$
 whose kernel consists of elements of $1+2 \End_{\Z_2}(T_2(X))$.  We have assumed that
the $G$-module $X_2$ is absolutely simple. This implies that the $H$-module $X_2$
is also absolutely simple. Here the structure of $H$-module is defined on
$X_2$ via
$$H\subset\Aut_{\Z_2}(T_2(X)) \twoheadrightarrow \Aut(X_2).$$
The absolute simplicity of the $H$-module $X_2$ means that the natural homomorphism
$$\F_2[H] \to \End_{\F_2}(X_2)$$
is surjective. By Nakayama's Lemma, this implies that the natural homomorphism
$$\Z_2[H] \to \End_{\Z_2}(T_2(X))$$
is also surjective (see \cite{Mazur}, p. 252).
 
Let $V_2(X)=T_2(X)\otimes_{\Z_2}\Q_2$ be the $\Q_2$-Tate module of $X$. It is
well-known that $V_2(X)$ is the $2g$-dimensional $\Q_2$-vector space and
$T_2(X)$ is a $\Z_2$-lattice in $V_2(X)$.
Clearly, the $\Q_2[H]$-module $V_2(X)$ is also absolutely simple.

The choice of polarization on $X$ gives rise to a non-degenerate alternating
bilinear form (Riemann form) \cite{MumfordAV}
$$e: V_{2}(X) \times V_2(X) \to \Q_2(1) \cong \Q_2.$$
Since $F$ contains all $2$-power roots of unity, $e$ is $\Gal(F)$-invariant
and therefore is $H$-invariant. This means that $H$ is a subgroup of the
corresponding symplectic group $\Sp(V_2(X),e)$. We have
$$H \subset \Sp(V_2(X),e) \cong \Sp_{2g}(\Q_2)\subset \Sp_{2g}(\bar{\Q}_2).$$

There exists a finite Galois extension $L$ of $K$ such that all
endomorphisms of $X$ are defined over $L$. We write $\End^0(X)$ for the
$\Q$-algebra $\End(X)\otimes\Q$ of endomorphisms of $X$.
Since $X$ is supersingular,
$$\dim_{\Q}\End^0(X)=(2\dim(X))^2=(2g)^2.$$
Recall (\cite{MumfordAV}) that the natural map
$$\End^0(X)\otimes_{\Q}\Q_{2} \to \End_{\Q_{2}}V_{2}(X)$$
is an embedding.
Dimension arguments imply that
$$\End^0(X)\otimes_{\Q}\Q_{2} = \End_{\Q_{2}}V_{2}(X).$$
Since all endomorphisms of $X$ are defined over $L$, the image
$$\rho_{2,X}(\Gal(L)) \subset \rho_{2,X}(\Gal(F)) \subset\Aut_{\Z_2}(T_2(X))
\subset\Aut_{\Q_2}(V_2(X))$$
commutes with $\End^0(X)$. This implies that
$\rho_{2,X}(\Gal(L))$ commutes with  $\End_{\Q_{2}}V_{2}(X)$ and therefore
consists of scalars. Since
$$\rho_{2,X}(\Gal(L)) \subset \rho_{2,X}(\Gal(F)) \subset \Sp(V_2(X),e),$$
$\rho_{2,X}(\Gal(L))$ is a finite group. Since $\Gal(L)$ is a subgroup of
finite index in $\Gal(F)$, the group $H=\rho_{2,X}(\Gal(F))$ is also finite.
In particular, the kernel of the reduction map modulo $2$
$$\Aut_{\Z_2}T_2(X) \supset H \twoheadrightarrow G \subset \Aut(X_2)$$
consists of  elements of finite order and, thanks to the Minkowski-Serre Lemma,
$Z:=\ker(H \to G)$ has exponent $1$ or $2$. In particular, $Z$ is
commutative. We have
$$Z\subset H \subset \Sp(V_2(X),e) \cong \Sp_{2g}(\Q_2)\subset \Sp_{2g}(\bar{\Q}_2).$$
Since $Z$ consists of semisimple elements and rank of $\Sp_{2g}$ is $g$,
the group $Z$ is isomorphic (``conjugate") to a multiplicative subgroup of $(\bar{\Q}_2^*)^g$. Since the exponent of $Z$ is either $1$ or $2$, the group
$Z$ is isomorphic to a multiplicative subgroup of $\{1, -1\}^g$. Hence $Z$ is an $\F_2$-vector space of dimension $d\le g$. This implies that the adjoint action
$$H \to H/Z=G \to \Aut(Z) \cong\GL_d(\F_2)$$
is trivial, since every nontrivial representation of $G$ in characteristic $2$ must have dimension strictly greater than $g \ge d$.
This means that $Z$ lies
in the center of $H$. Since the $\Q_2[H]$-module $V_2(X)$ is faithful and absolutely
simple, $Z$ consists of scalars. This implies that either $Z=\{1\}$ or
$Z=\{\pm 1\}$. In other words, either $H\cong G$ or $H \twoheadrightarrow G$ is a double cover. In both cases $V_2(X)$ is an absolutely irreducible representation of $H$ of dimension $2g$ over $\Q_2$. But by Proposition \ref{double} applied to $G'=H$ and $V=V_2(X)$,
$$\dim_{\Q_2}(V_2(X)) \ne 2g.$$
This gives us the desired contradiction. This ends the proof of Lemma \ref{supernot} and therefore of Theorems \ref{main2} and \ref{M12}.
\end{proof}

\begin{ex}
Suppose $p$ is an odd prime, $q>1$ is a power of $p$, $m>2$ is an even integer. Let us put $n=(q^m-1)/(q-1)$.
Suppose $k$ is an algebraically closed field of characteristic $p$ and $K=k(z)$ is the field of rational functions. The Galois group of $x^m+zx+1$ over $K$ is $\L_m(q)$ and the Galois group of  $x^m+x+z$ over $K$ is $\PGL_m(\F_q)$ (\cite{Ab}, p. 1643). Therefore the jacobians of the hyperelliptic curves
$y^2=x^m+zx+1$ and $y^2=x^m+x+z$ have no nontrivial endomorphisms over an algebraic closure of $K$.
\end{ex}


\begin{thebibliography}{99}

\bibitem{Ab} S. S. Abhyankar, {\em Projective polynomials}. Proc. AMS {\bf 125} (1997), 1643--1650.

\bibitem{Atlas} J. H. Conway, R. T. Curtis, S. P. Norton, R. A. Parker, R. A. Wilson, Atlas of finite groups. Clarendon Press, Oxford, 1985.

\bibitem{GurTiep} R. M. Guralnick, Pham Huu Tiep, {\em Low-dimensional representations of special linear groups in cross characteristic}. Proc. London Math. Soc. {\bf 78} (1999), 116--138.

\bibitem{Hump} J. F. Humphreys, {\em The projective characters of the Mathieu group} $\M_{12}$ {\em and of its automorphism group}. Math. Proc. Camb. Phil. Soc. {\bf 87} (1980), 401--412,

\bibitem{James} G. D. James, {\em The modular characters of the Mathieu 
groups}.    J. Algebra {\bf 27} (1973), 57--111.

\bibitem{AtlasB} Ch. Jansen, K. Lux, R. Parker, R. Wilson,  An Atlas
of Brauer characters. Clarendon Press, Oxford, 1995.


\bibitem{Katz1} N. Katz, {\em Monodromy of families of curves:
    applications of some results of Davenport-Lewis}. In:
    S\'eminaire de Th\'eorie des Nombres, Paris 1979-80
    (ed. M.-J. Bertin); Progress in Math. {\bf 12},
    pp. 171--195, Birkh\"auser, Boston-Basel-Stuttgart,
    1981.

\bibitem{Katz2} N. Katz,  {\em Affine cohomological transforms,
    perversity, and monodromy}.
     J. Amer. Math. Soc. {\bf 6} (1993), 149--222.

\bibitem{Masser} D. Masser, {\em Specialization of
    some hyperelliptic jacobians}. In:
    Number Theory in Progress
    (eds.  K. Gy\"ory, H. Iwaniec, J. Urbanowicz), vol. I, pp. 293--307;
     de Gruyter, Berlin-New York, 1999.

\bibitem{Mazur} B. Mazur, {\em Deformation theory of Galois representations}. In: Modular forms and Fermat's last theorem (G. Cornell, J. H. Silverman, G. Stevens, eds.), pp. 243--311, Springer-Verlag, New York, 1997.


\bibitem{MumfordAV} D. Mumford, Abelian varieties, Second edition,
 Oxford University Press, London, 1974.

\bibitem{Mori1} Sh. Mori, {\em The endomorphism rings of some abelian varieties}.  Japanese J. Math,  {\bf 2}(1976), 109--130.

\bibitem{Mori2} Sh. Mori, {\em The endomorphism rings of some abelian varieties}. II, Japanese J. Math,  {\bf 3}(1977), 105--109.



\bibitem{Zarhin} Yu. G. Zarhin, {\em Hyperelliptic jacobians without
complex multiplication}. Math. Res. Letters {\bf 7}(2000), 123--132.
\bibitem{Zarhin2} Yu. G. Zarhin, {\em Hyperelliptic jacobians and modular representations}, http://xxx.lanl.gov/abs/math.AG/0003002, to appear in
Texel volume ``Moduli of abelian varieties" (eds. G. van der Geer, C. Faber, F. Oort), Birkh\"auser.


\end{thebibliography}
\end{document}